\documentclass{amsart}
\usepackage{graphicx} 
\usepackage{amsmath,amssymb, amsfonts, amscd}
\usepackage{tikz-cd, tikz}
\usepackage{subcaption}
\usepackage{amsthm, bbm}
\usepackage{appendix}
\usepackage{enumitem}
\usepackage[normalem]{ulem}
\usepackage{hyperref}
\usepackage[capitalise]{cleveref}
\usepackage[a4paper, margin=3cm]{geometry}
\usepackage{xfrac}
\usepackage{mathtools,tensor}
\usepackage[style=alphabetic, backend=biber]{biblatex}
\usepackage{xcolor}
\usepackage{framed}

\theoremstyle{plain}
\newtheorem*{theorem*}{Theorem}
\newtheorem*{remark*}{Remark}
\newtheorem*{example*}{Example}
\newtheorem*{conjecture*}{Conjecture}

\newtheorem{lemma}{Lemma}[subsection]
\newtheorem{corollary}[lemma]{Corollary}
\newtheorem{theorem}[lemma]{Theorem}

\newtheorem{thm}[lemma]{Theorem}

\newtheorem{prop}[lemma]{Proposition}

\newtheorem{crl}[lemma]{Corollary}
\newtheorem{conj}[lemma]{Conjecture}

\newtheorem{introtheorem}{Theorem}

\newtheorem{introthm}[introtheorem]{Theorem}

\newtheorem*{introquestion*}{Question}
\newtheorem*{introremark*}{Remark}
\newtheorem{definition}[lemma]{Definition}
\newtheorem*{definition*}{Definition}

\newtheorem{example}[lemma]{Example}

\theoremstyle{remark}
\newtheorem{remark}[lemma]{Remark}

\newtheorem{notation}[lemma]{Notation}
\newtheorem{rem}[lemma]{Remark}

\newcommand{\C}{\mathcal{C}}
\newcommand{\1}{\mathbbm{1}}
\newcommand{\D}{\mathcal{D}}
\newcommand{\E}{\mathcal{E}}
\newcommand{\M}{\mathcal{M}}
\newcommand{\N}{\mathbb{N}}
\newcommand{\Z}{\mathbb{Z}}
\newcommand{\End}{\mathrm{End}}
\newcommand{\FPdim}{\mathrm{FPdim}}
\newcommand{\id}{\mathrm{id}}
\newcommand{\Rep}{\mathrm{Rep}}
\newcommand{\Vect}{\mathrm{Vec}}
\newcommand{\Hom}{\mathrm{Hom}}
\newcommand{\kk}{\mathbbm{k}}
\newcommand{\pic}{\mathrm{Pic}}

\newcommand{\eps}{\epsilon}

\newcommand{\res}{\mathrm{Res}}

\newcommand{\comment}[1]{}
\newcommand{\prescr}[2]{\tensor*[^{#1}]{#2}{}}

\title{Generalized near-group tensor categories}
\author{Inna Entova-Aizenbud, Daniel Sebbag}
\date{\today}

\begin{document}

\begin{abstract}
   Let $\D$ be a finite tensor category with exactly one isomorphism class of simple objects which are non-invertible. 
   
   Semisimple tensor categories $\D$ with this property are called near-group categories, and have been studied by Evans, Gannon, Siehler and others. In this paper, we study such categories in the non-semisimple setting, with the additional assumption that the non-invertible simple object is projective. Specifically, we consider the case when $\D$ admits a fiber functor and has a unique semisimple abelian block (in the general case that $\D$ is not semisimple). We describe the structure of such categories and their corresponding Hopf algebras, and construct examples of such categories with a cyclic Picard group. We obtain additional results when the Frobenius-Perron dimension of the unique non-invertible simple object is a prime number.
\end{abstract}

\maketitle

\section{Introduction}
\subsection{}
    In this paper, we consider a special class of finite tensor categories: those with exactly one simple object which is non-invertible. We will call such categories {\it generalized near-group categories}.

    To understand such categories, it makes sense to consider the following data associated with $\D$:
    \begin{enumerate}
        \item The Picard group $G:=\pic(\D)$ of $\D$: this is the group whose elements are isomorphism classes of invertible objects. The group multiplication corresponds to tensor products of invertible objects.
        \item Let $Q\in \D$ be the non-invertible simple  object, and denote by $r$ the multiplicity of the simple object $Q$ among the composition factors of $Q\otimes Q$. We will call $r$ the {\it extension constant} of $\D$.
    \end{enumerate}
    It is not difficult to see that when $Q$ is also projective, the Grothendieck ring of $\D$ is completely determined by $G$, $r$ and $\FPdim(P(\1))$, where $P(\1)$ is the projective cover of the unit object $\1$. One may pose the following question: 
    \begin{introquestion*}\label{ques:main}
        Consider a triplet $(G,r,d)$, where $G$ is a finite group, $r,d\in \Z_{\geq 0}$.
        \begin{enumerate}
            \item For which triplets $(G,r,d)$ can one find a finite tensor category $\D$ associated with this triplet? 
            \item Is there a finite number of equivalence classes of finite tensor categories associated with a given triplet $(G,r,d)$?
        \end{enumerate}
    \end{introquestion*}
The pair $(G, r)$ is called the {\it generalized fusion rule} of $\D$.

    The generalized near-group categories which are semisimple (fusion categories) have been studied in the literature under the name near-group categories. When $r=0$, these are also known as Tambara-Yamagami categories. They were first defined in \cite{TY}, generalizing the Ising model, and studied further in \cite{Si1, Si2, EG2, T2, T3, O, Iz} and more\footnote{In \cite{T3}, the term \textit{generalized near-group category} refers to a fusion category so that the Picard group acts transitively on the set of isomorphism classes of non-invertible simple objects.}; a list of many such categories and a partial classification appear in \cite{EG2}. Tannakian near group categories over $\mathbb{C}$ are precisely the categories of finite-dimensional representations of finite groups which have a unique irreducible representation of dimension larger than $1$. This class of groups was classified in characteristic zero by \cite{Sei}.

    When working over an algebraically closed field of characteristic $0$, there are two types of finite groups with a unique irreducible representation of dimension greater than $1$ (see \cite{Sei}):
    \begin{itemize}
        \item An extra-special group of order $2^{1+2n}$ for $n \geq  1$ (there are two of each order), which is a central extension of $(\Z/2\Z)^{2n}$ by $\Z/2\Z$.
        \item The finite affine group $Aff_1(\mathbb{F}_q)=\mathbb{F}_q\rtimes \mathbb{F}_q^{\times}$.
    \end{itemize}
 
    These two families are also ``typical" examples for near-group categories, and the unique non-invertible simple object imposes strict limitations on the structure of the categories. In particular, the combined results of \cite{EG2, Iz, TY} show that for a unitary near-group category with fusion rule $(G,r)$, the group $G$ must be either abelian or an extra-special $2$-group, and if $r>|G|$ then $|G|\rvert r$.
    
    As a starting point in the study of generalized near-group categories, we are interested in generalized near-group categories possessing two additional properties:
    \begin{itemize}
        \item $\D$ possesses a fiber functor into $\Vect$, making $\D$ equivalent to the category of finite-dimensional modules over some finite-dimensional Hopf algebra $H$. In particular, in this case $\D$ is integral. We assume that this fiber functor is \textit{flip compatible} (see \cref{def: flip-compatible fiber functor}).
        \item $\D$ is {\it almost-semisimple}: that is, the non-invertible simple object $Q$ is projective (consequently, injective).
    \end{itemize}
    An important feature of almost-semisimple generalized near-group categories is that the Serre subcategory generated by the invertible objects of $\D$ is a direct summand of $\D$. 

    In this paper, we focus on the case when $char(\kk) = 0$. Our main result is:
    \begin{introthm}\label{main theorem}
        Let $\D$ be an almost semisimple generalized near-group category admitting a flip compatible fiber functor $F:\D\to \Vect$, and let $\C\subset \D$ be the Serre subcategory generated by the invertible objects.

        Let $H:=\End(F)$. We have an algebra isomorphism $$H \cong A\#_\sigma K$$ where $A\subset H$ denotes the smallest non-trivial normal Hopf subalgebra of $H$, while $K := \End(F|_{\C})$ and $\sigma:K\otimes K\to A$ is a $2$-cocycle. Moreover, the algebra $A$ is semisimple.

        Equivalently, there exists a canonical exact sequence of tensor categories:
        $$
            \C\to \D\to \E
        $$
        where
        $\E=\Rep(A)$ is a fusion category.
    \end{introthm}
    \begin{introremark*}
        A priori, the category $\E$ depends on the fiber functor $F$, but $\FPdim(\E)$ depends only on $\FPdim(\D)$ and $\FPdim(\C)$.    
    \end{introremark*}
    
    From \cref{main theorem} we get:
    \begin{introthm}        
        If $\D$ is an almost semisimple generalized near-group category admitting a fiber functor with $r=0$, then $A = \kk[\Z/2\Z]$ is central in $H$.
    \end{introthm}
    \begin{introremark*}
        When $r=0$, the flip-compatibility of the fiber functor is automatic.
    \end{introremark*}
    In \cref{sec: Extending H_{n,d}}, we construct a family of almost semisimple generalized near-group categories admitting a fiber functor: for every $n\geq 2$, we construct such a category with $G=\Z/n\Z$, $r=0$. This is done using certain extensions of Taft algebras. 
    
    Our construction demonstrates the added complexity of the non-semisimple case: in the fusion (semisimple) case, for an integral near-group category with fusion rule $(G,0)$, the size $|G|$ must be a perfect square (even without having a fiber functor). The example we construct shows that the non-semisimple case is more diverse.

    In the categorical point of view, an almost semisimple generalized near-group category is a special case of the following construction: given a finite tensor category $\C$ and an indecomposable exact $\C$-bimodule category $\M$, one may ask whether it is possible to define a tensor structure on $\C\oplus \M$ that is compatible with the tensor structure of $\C$ and the $\C$-bimodule structure of $\M$. \cref{main theorem} suggests that the answer to this question (at least in some cases) may be resolved using exact sequences of tensor categories (as defined by \cite{BN} and extended by \cite{EG1} for cases when there is no fiber functor).

    A future direction of research would be understanding almost semisimple generalized near-group categories in positive characteristic. The first natural example to look at would be Tannakian categories (i.e. finite tensor categories of the form $\Rep_{\kk}(G)$ where $G$ is a finite group scheme and $\kk$ is a field of characteristic $p$). 
    
    The classification of finite groups $G$ for which $\Rep_{\kk}(G)$ is a generalized near-group category follows from the results of \cite{DN}. In \cite{DN}, Dolfi and Navarro prove that for a group $G$ such that the largest normal $p$-subgroup is trivial, $G$ has a unique irreducible representation of dimension greater than $1$ if and only if $G$ is one of the following groups (in particular, $G$ is solvable):
    \begin{enumerate}
        \item $G$ is a doubly transitive Frobenius group, with Frobenius complement $H$ such that $H'$ is a $p$-group;
        \item $p\ne 2$ and $G$ is an extra-special $2$-group;
        \item $p=2$ and $G = a\Gamma(q^2)$ is the affine semilinear group, with $q$ a Mersenne prime;
        \item  There are two more groups arising when $p=2$. 
    \end{enumerate}
   
We believe that over a field of positive characteristic, even the classification of generalized near group categories of small Frobenius-Perron dimension could be a very interesting and non-trivial problem.
\subsection{Structure of the paper}

   \cref{sec: Preliminaries} contains the preliminaries.
    In Sections \ref{sec:general}, \ref{sec: integral near-group categories}, we describe some basic results concerning almost-semisimple near-group categories, both with and without assuming the existence of a fiber functor. In \cref{sec: integral near-group categories}, we prove \cref{main theorem}. 
    
    In \cref{sec: Extending H_{n,d}}, we construct explicit examples of almost semisimple generalized near-group categories admitting a fiber functor, based on Taft algebras. 

    In \cref{sec: classification first results}, we study the case when $\FPdim(Q)$ is a prime number.

\subsection{Acknowledgements}
We thank David Benson for useful discussions and for pointing out to us the results \cite{DN}. 

This research was supported by the Israel Science Foundation grant ISF No. 1362/23 (PI: Inna Entova-Aizenbud). 
Claude Opus 5.5 was used to verify a late draft of this paper. The final version incorporates corrections and improved exposition based on the AI report. 

\section{Preliminaries}\label{sec: Preliminaries}
    Throughout this paper, we work over $\kk$ an algebraically closed field of characteristic zero, and all categories and functors are assumed to be linear over $\kk$. All the Hopf algebras will be finite-dimensional, as will be their modules and comodules.

    \subsection{Notation and terminology}\label{ssec:notn}
        A category $\D$ is called a {\it finite tensor category} if $\D$ is a rigid monoidal abelian $\kk$-linear category with the following properties:
        \begin{itemize}
            \item the bifunctor $-\otimes -$ is bilinear,
            \item all $\Hom$-spaces in $\D$ are finite dimensional and all objects have finite length,
            \item $\D$ has a finite number of (isomorphism classes of) simple objects,
            \item $\D$ has enough projective objects,
            \item the unit object $\1$ is simple.
        \end{itemize} Throughout this paper, we assume that all categories have strict associativity and unit constraints. 

        Under these assumptions, every projective object is also injective (see \cite[Proposition 6.1.3]{EGNO}), and the tensor product is a biexact functor.

        For a finite tensor category $\D$, a \textit{fiber functor} on $\D$ is an exact monoidal functor $F:\D\to \Vect$. 
        A fiber functor $\D\to \Vect$ defines equivalences of tensor categories $\D\xrightarrow{\sim} H-Comod\xrightarrow{\sim} H^*-mod$ into the category of finite-dimensional left $H$-comodules (equivalently into the category of finite-dimensional left $H^*$-modules) where $H:=Coend(F)$ is a finite-dimensional Hopf algebra and $H^*\cong \End(F)$ is its dual (see \cite[Theorem 5.3.12]{EGNO}).

        An object $X\in{\D}$ is called {\it invertible} if there exists a simple object $Y\in \D$ such that $X\otimes Y = \1$. The Picard group of $\D$, denoted by $Pic(\D)$, is defined to be the set of isomorphism classes of invertible objects, with multiplication operation given by the tensor product.

        A finite tensor category $\D$ is called \textit{pointed} if every simple object $X\in{\D}$ is invertible. Similarly, a Hopf algebra $H$ is called pointed if the category of finite-dimensional $H$-comodules is pointed.

        For finite tensor categories $\C,\D$, a tensor functor $F:\C\to \D$ is a monoidal exact functor. Such a functor is called \textit{surjective} if every object $Y\in \D$ is a subquotient of $F(X)$ for some object $X\in \C$.
        
        For a Hopf algebra $H$, we denote by $H^+$ its augmentation ideal, i.e. $H^+:=\ker(\eps)$. For a Hopf algebra $H$ and $A<H$ a Hopf subalgebra, $A$ is called \textit{normal} if $A$ is stable under the left and right adjoint action of $H$. 
    \subsection{Reconstructing Hopf algebras}
        Let $\D$ be a finite tensor category with a fiber functor $F: \D\to \Vect$. Then $H:=Coend(F)$ has a natural structure of a (finite-dimensional) Hopf algebra and $\D\cong H-Comod$ (see \cite[Sections 1.10, 5.3]{EGNO}), with multiplication induced by tensor product. Explicitly, $H=\left(\bigoplus_{X\in \D} F(X)^*\otimes F(X)\right)/I$ where $$I = span\left\{~y^*\otimes F(f)(x)-F(f)^*(y^*)\otimes x ~~| ~~X, Y\in \D,~ f\in \hom_{\D}(X,Y),~ x\in F(X), ~y^*\in F(Y)^*~\right\}.$$
 
        For each object $X\in \D$ we have the canonical coaction of $H$ on $F(X)$. This gives a coalgebra homomorphism
        $\rho_X:F(X)\otimes F(X)^* \to H$, whose image will be denoted by $E(X)$. We have the following properties of the coalgebras $E(X)$:
        \begin{lemma}\label{lem:E_X}
            \mbox{}
            \begin{enumerate}           
                \item The coalgebra $H$ is a direct limit of coalgebras $ E(X)$ for $X\in \D$. 
                \item For any $X, Y\in \D$ we have: $E(X\oplus Y)=E(X)+E(Y)$. In particular, for any $m\in \Z_{\geq 0}$, we have: $E(X^{\oplus m})=E(X)$.
                \item If $f:X\to Y$ is surjective, then $E(Y)\subset E(X)$.
                \item If $X$ is a simple object in $\D$ then $E(X)$ is a simple coalgebra, and $E(X)\cong F(X)\otimes F(X)^*$ as coalgebras.
                \item If $P_1, \ldots, P_k$ denote the isomorphism classes of indecomposable projective objects in $\D$, then $H=\sum_i E(P_i)$.           
                \item The multiplication operation $H\otimes H\to H$ sends $E(X)\otimes E(Y)\to E(X\otimes Y)$.
            \end{enumerate}
        \end{lemma}
        \begin{proof}
            Statements (1)-(4) follow immediately from the explicit description of $H$ given above, and statement (6) follows from the definition of multiplication on $H$. Finally, (5) follows directly from (1), (2), (3), since every object in $\D$ has a projective cover which is, in turn, a direct sum of indecomposable projective objects.
        \end{proof}

    \subsection{Left and right coinvariants of Hopf algebra homomorphisms}\label{sec: left and right coinvariances}
        In this subsection we recall the notions related to short exact sequences of Hopf algebras \cite[\S1.7]{M}.
\begin{definition}[Coinvariants]
    Let $\pi:H\to K$ be a homomorphism of Hopf algebras. The sets of {\it right and left $\pi$-coinvariants} are defined by $$H^{co \pi}:=\{h\in H~\rvert~ (\id\otimes \pi)\Delta(h) = h\otimes 1\}, \;\;\;\;\;\prescr{co\pi}{H} := \{h\in H~\rvert ~(\pi\otimes \id)\Delta(h) = 1\otimes h\}.$$ 
\end{definition}
It is easy to see that $H^{co \pi}, \prescr{co\pi}{H}$ are subalgebras of $H$, and that $H^{co \pi}$ is a left coideal, while $\prescr{co\pi}{H}$ is a right coideal. Moreover, the antipode $S$ gives a bijection between $H^{co \pi}$ and $ \prescr{co\pi}{H}$. 

\begin{lemma}\label{lem:coinvar_kernel}
Let $\pi:H\to K$ be a homomorphism of Hopf algebras. Then we have:
\begin{enumerate}
    \item $\pi(h)=\varepsilon(h)1$ for any $h\in H^{co \pi}\cup \prescr{co\pi}{H}$. In particular, $ (H^{co \pi})^+, (\prescr{co\pi}{H})^+\subset  \ker \pi.$
    \item Let $A\subset H$ be a sub-bialgebra such that $\pi(h)=\varepsilon(h)1$ for any $h\in A$. Then $A\subset H^{co \pi}\cap \prescr{co\pi}{H} $.
\end{enumerate}
\end{lemma}
\begin{proof}
\mbox{}
\begin{enumerate}
    \item 
We prove the statement for $h\in H^{co \pi}$ (the statement for $h\in \prescr{co\pi}{H}$ is proved analogously, or by applying the antipode map to $h$).
For $h\in H^{co \pi}$, we have $$\Delta(\pi(h))=(\pi\otimes\id)(\mathrm{id}\otimes\pi)\Delta (h)=(\pi\otimes\id)(h\otimes 1)=\pi(h)\otimes 1.$$ Applying $\varepsilon\otimes\id$ to both sides yields $\pi(h)=\varepsilon(h)1$ as required.
    \item Let $h\in A$. Then $\Delta(h)\in A\otimes A$ so $$(\id\otimes \pi)\Delta(h) = (\id\otimes \varepsilon(-)1)\Delta (h) = h\otimes 1, \;\;\;\;\;(\pi\otimes \id)\Delta(h) = ( \varepsilon(-)1\otimes \id)\Delta (h)= 1\otimes h,$$
so $A\subset H^{co \pi}\cap \prescr{co\pi}{H} $.
\end{enumerate}
\end{proof}
Thus both left and right coinvariants play the role analogous to the kernel of a group homomorphism. Under a suitable condition, we may make this more precise.

In fact, it is possible to show that when $\dim H, \dim K<\infty$, the $2$-sided ideal $\ker \pi$ is generated, both as a left ideal and a right ideal, by each of the spaces $ (H^{co \pi})^+, (\prescr{co\pi}{H})^+$.

Finally, we mention the following corollary of the Nichols–Zoeller theorem (\cite{NZ}), applied to the embedding $K^*\hookrightarrow H^*$:
\begin{lemma}[See \cite{Sch}]\label{lem:dim_coinvar}
   Let $\pi:H\to K$ be a surjective homomorphism of finite-dimensional Hopf algebras. Then $\dim H^{co \pi}=\dim \prescr{co\pi}{H}=\dim H/\dim K$. 
\end{lemma}
 
\subsection{Coinvariants and short exact sequences of Hopf algebras}
\begin{definition}
    For Hopf algebras $A,H,K$, a short exact sequence of Hopf algebras
    $$
        \kk\to A\xrightarrow{i}H\xrightarrow{\pi} K\to \kk
    $$
    is a sequence in which $i$ is injective, $\pi$ is surjective, $i(A) = H^{co \pi}$ and $\pi\circ i = \eps$, $\ker \pi = Hi(A)^+$.
\end{definition}

\begin{lemma}\label{lem:normal_equiv_cond}
    Let $\pi:H\to K$ be a surjective homomorphism between finite-dimensional Hopf algebras. The following conditions are equivalent:
    \begin{enumerate}
        \item $H^{co \pi} = \prescr{co\pi}{H}$. 
        \item $H^{co \pi}$ is a normal Hopf subalgebra of $H$. 
        \item The antipode $S$ preserves the subspace $H^{co \pi}$. 
       \item There exists a Hopf subalgebra $A\subseteq H$ with $\pi|_A(\cdot)=\varepsilon(\cdot)1$ and  $\dim A\cdot\dim K=\dim H$.
        \item There exists a Hopf subalgebra $A\subseteq H$ with $\pi|_A(\cdot)=\varepsilon(\cdot)1$ and  $\ker \pi = HA^+$.
        \item The following is an exact sequence of Hopf algebras:
        $$
            \kk \to H^{co \pi}\to H\xrightarrow{\pi} K\to \kk.
        $$ 
      
    \end{enumerate}    
\end{lemma}
\begin{proof}
We first prove that $(1)\Leftrightarrow(3)\Leftrightarrow(2)$ and then $(2)\Rightarrow (4)\Leftrightarrow (5)\Rightarrow (6) \Rightarrow(2)$.

\begin{enumerate}

    \item[\textbf{[(1) $\Longleftrightarrow$ (3)]}] Follows from the fact that the antipode $S$ gives a bijection between $H^{co\pi}$ and $\prescr{co\pi}{H}$.
    
     \item[\textbf{[(1) $\Longrightarrow$ (2)]}] Assume $H^{co \pi} = \prescr{co\pi}{H}$.    
     Recall that $H^{co\pi}$ is a left coideal subalgebra and $\prescr{co\pi}{H}$ is a right coideal subalgebra. When these subalgebras coincide, $H^{co\pi}$ is both a left and a right coideal, so it is a subcoalgebra of $H$. Since $H^{co\pi}$ is preserved by $S$ due to (3), we conclude that $H^{co\pi}$ is a Hopf subalgebra of $H$.

     It remains to prove normality.
 Let $a\in H^{co \pi}$, $h\in H$. Then $$(\mathrm{id}\otimes\pi)\Delta\big(\sum h_1aS(h_2)\big)=\textstyle\sum h_1aS(h_4)\otimes\pi(h_2S(h_3))=(\sum h_1aS(h_2))\otimes1$$
 where $\Delta (h) =: \sum h_1\otimes h_2$ and $(\Delta\otimes \id^{\otimes 2})(\Delta\otimes \id)\Delta(h)=:\sum h_1\otimes h_2\otimes h_3 \otimes h_4$.
 So $\sum h_1aS(h_2)\in H^{co \pi}$ and thus $H^{co \pi}$ is stable under the left adjoint action of $H$. The stability under the right adjoint action is proved analogously.

     \item[\textbf{[(2) $\Longrightarrow$ (3)]}]
If $H^{co\pi}$ is a Hopf subalgebra of $H$, then is it preserved by the antipode $S$.

\item[\textbf{[(2) $\Longrightarrow$ (4)]}] This implication follows directly from \cref{lem:coinvar_kernel} and \cref{lem:dim_coinvar}.
     \item[\textbf{[(4) $\Longleftrightarrow$ (5)]}] 
      
Now, by our assumption on $A$, we have:  $HA^+\subset \ker\pi$.
Let us show that $ HA^+= \ker\pi$ if and only if $\dim A \cdot \dim K=\dim H$.

By the Nichols-Zoeller Theorem (\cite{NZ}), $H$ is free over $A$. Let $m:=  rk_A H=\frac{\dim(H)}{\dim(A)} \in \Z$.

Consider the short exact sequence of left $A$-modules
        $$
            0\to A^+\to A\xrightarrow{\varepsilon} \1_A\to 0.
        $$ where $\1_A $ stands for the trivial $1$-dimensional $A$-module.
        Applying $H\otimes_A-$ gives a short exact sequence of $H$-modules
        $$
            0\to H A^+\to H\to H\otimes_A \1_A\to 0
        $$
        and $\dim H\otimes_A \1_A = rk_A H=m$. Since $ HA^+\subset \ker\pi$, we obtain a surjective map of left $H$-modules $ H\otimes_A \1_A\twoheadrightarrow K$.  
     We conclude that $$ m = \frac{\dim(H)}{\dim(A)}=\dim H-\dim HA^+ \geq  \dim K$$ and there is an equality here if and only if $HA^+= \ker\pi $.    
      \item[\textbf{[(5) $\Longrightarrow$ (6)]}] 
      Assume (5) holds. By the above argument, $\dim A \cdot \dim K=\dim H$. \cref{lem:dim_coinvar} now implies that $A= H^{co \pi} =\prescr{co\pi}{H}$, so (6) holds. 
     
 \item[\textbf{[(6) $\Longrightarrow$ (2)]}] By \cref{lem:coinvar_kernel} and \cref{lem:dim_coinvar} (6)$\Longrightarrow$(1) and we proved that (1)$\Leftrightarrow$(2).
      
  
\end{enumerate}
\end{proof}

\begin{definition}
      A surjective homomorphism of Hopf algebras $\pi:H\to K$ is called {\it normal} if the conditions of \cref{lem:normal_equiv_cond} hold.
\end{definition}
\cref{lem:coinvar_kernel} and \cref{lem:normal_equiv_cond} now imply:
\begin{corollary}
     Let $\pi:H\to K$ be a normal homomorphism between finite-dimensional Hopf algebras. Let $M$ be an $H$-module. The action of $H$ on $M$ factors through $\pi$ if and only if $M$ is a trivial $H^{co \pi}$-module (i.e. any $h\in H^{co \pi}$ acts by $\eps(h)\id_M$). 
\end{corollary}

  \subsection{Crossed product of Hopf algebras}\label{ssec:crossed_prod}
        \begin{definition}
            For a pair of Hopf algebras $(A,K)$, we say that $K$ {\it measures} $A$ if there exists a linear map $K\otimes A\to A$, $k\otimes a\mapsto k. a$ so that the following conditions hold:
        \begin{itemize}
            \item $k. 1_A = \eps_K(k)1_A$ for all $k\in K$, where $\eps_K$ is the counit in $K$.
            \item $k. (ab) = \sum (k_1. a)(k_2. b)$ for all $k\in K$, $a,b\in A$, where $\Delta_K (k) = \sum k_1\otimes k_2$ is the comultiplication in $K$.
        \end{itemize}
        \end{definition}

        Fix a pair of Hopf algebras $(A,K)$ such that $K$ measures $A$, and a $2$-cocycle $\sigma: K\otimes K\to A$. We may assume $\sigma$ is {\it normalized}, meaning that $\sigma(h,1_K) = \sigma(1_K,h) = \eps(h)1_A$ for all $h\in K$.
        \begin{definition}
             The {\it crossed product} of $A$ with $K$, denoted by $A\#_\sigma K$, is the  vector space $A\otimes K$ with multiplication defined by:
        $$
            (a\otimes h)\cdot (b\otimes k) = \sum a(h_1.b)\sigma(h_2, k_1)\otimes h_3k_2 \;\; \text{ for }\;\; a,b\in A, \;h,k\in K.
        $$
        where $\Delta_K (k) = \sum k_1\otimes k_2$ and $(\Delta_K\otimes \id)\Delta_K(h)=\sum h_1\otimes h_2\otimes h_3 $.
        \end{definition}
        We denote the elements of $A\#_\sigma K$ by $a\#k$ instead of $a\otimes k$. 
        
        Now, assume that $A$ is a twisted $K$-module, i.e. $1_K. a = a$ for all $a\in A$ and 
        $$
            h.(k.a) = \sum \sigma(h_1,k_1)(h_2k_2.a)\sigma^{-1}(h_3,k_3)
        $$
        $\forall h,k\in K$,  $a\in A$. Then, $A\#_\sigma K$ becomes an associative algebra with unit $1\#1$ \cite[Lemma 7.1.2]{M}.
       
        If $\pi:H\to K$ is a normal surjective homomorphism of finite dimensional Hopf algebras, then by \cite{Sch}, $H$ is isomorphic to the crossed product of $H^{co \pi}$ and $K$. That is, there exists a $2$-cocycle $\sigma:K\otimes K\to H^{co \pi}$ such that $$H=H^{co \pi}\#_{\sigma}K.$$  
        This is proved using the notion of cleft extensions (\cite[Definition 7.2.1]{M}): using a convolution invertible $K$-comodule map $\gamma:K\to H$ (which is a section of $\pi$), one can construct a $2$-cocycle $\sigma$. For details, see \cite[\S7.2]{M}.
        
        \subsection{Exact sequences of tensor categories}\label{ssec:exact_seq_tensor_cats} In this subsection, we recall the notion of exact sequences of finite tensor categories, which was introduced in \cite{BN}. 

        Let $F:\C\to \D$ be a tensor functor between finite tensor categories. We will say that an object $Y\in \D$ is \textit{of $\Vect$-type} if it is a direct sum of finitely many copies of $\1_\D$. The functor $F$ is called \textit{normal} if any object $X\in \C$ admits a subobject $X'\in \C$ such that $F(X')$ is the largest $\Vect$-type subobject of $F(X)$. 
        
        We define $\mathfrak{Ker}_F$ to be the full subcategory of $\C$ whose objects $X$ satisfy: $F(X)$ is of $\Vect$-type in $\D$. The functor $F$ then defines a fiber functor $\mathfrak{Ker}_F\to \Vect$.

        Given three finite tensor categories $\C_1,\C_2,\C_3$ and tensor functors $i:\C_1\to \C_2$, $F:\C_2\to \C_3$,
        the sequence of tensor categories:
        $$
            \C_1\xrightarrow{i} \C_2 \xrightarrow{F} \C_3
        $$
        is called an \textit{exact sequence of tensor categories} if $F$ is surjective and normal, and $i$ induces an equivalence between $\C_1$ and $\mathfrak{Ker}_F\subset \C_2$. 
        \begin{remark}
            Notice that this definition requires in particular that $\C_1$ admits a fiber functor. The notion of exact sequences was generalized to avoid this restriction by \cite{EG1}, but in this work, we will consider exact sequences as defined in \cite{BN}.
        \end{remark}
        Notice that a short exact sequence of Hopf algebras
        $$
            \kk\to A \xrightarrow{i} H \xrightarrow{\pi} K\to \kk
        $$
        translates to an exact sequence of finite tensor categories 
        $$
            \Rep(K) \xrightarrow{F_1} \Rep(H)\xrightarrow{F_2} \Rep(A)
        $$
        for $F_1$ the pullback along $\pi$ and $F_2:=\res_A^H$. 
        \begin{itemize}
            \item The surjectivity of $\pi$ implies that $F_1$ is a full embedding.
            \item The injectivity of $i$ implies that $F_2$ is a surjective monoidal functor, in the sense of \cref{ssec:notn} (recall that $H$ is free over $A$ by Nichols-Zoeller theorem).

            \item Since $A\subset H$ is a normal Hopf subalgebra, the $A$-invariants subspace of any $V\in \Rep(H)$, given by $$V^A:=\{v\in V~\rvert ~a.v = \eps(a)v~\text{ for all }~a\in A\}$$ is an $H$-submodule. Thus $F_2$ is a normal functor.
            \item Recall that $\ker \pi = Hi(A)^+$. For every $V\in \Rep(H)$ such that $\res_A^H(V)$ is of $\Vect$-type, we have: $HA^+ V = 0$. So $V$ is an $H/HA^+$-module, and since $H/HA^+\cong K$, we conclude that $F_1$ is an equivalence between $\Rep(K)$ and $\mathfrak{Ker}_{F_2}$.
        \end{itemize}
   
\section{The basics of almost-semisimple generalized near-group categories}\label{sec:general}
    Let $\D$ be a generalized near-group category, and let $\C\subset \D$ be the (full) Serre subcategory of $\D$ generated by the invertible objects of $\D$. Let $Q$ be the (unique up to isomorphism) non-invertible simple object in $\D$.
    
The following lemma is straightforward:
    \begin{lemma}\label{lem:blocks_almost_semsimple}
        We have:
        \begin{enumerate}
            \item The subcategory $\C$ is a tensor subcategory of $\D$.
            \item If $\D$ is almost-semisimple, then $\C$ is a direct summand of $\D$ in the following sense: for any $M\in \D$, there exists a unique (up to isomorphism) $M'\in \C$ such that $M\cong M'\oplus Q^{\oplus m}$ for some $m\in \Z_{\geq 0}$. 
        \end{enumerate}
    \end{lemma}
       Assume now that $\D$ is almost-semisimple. 
       
       We have an equivalence of abelian categories $\D=\C\oplus \langle Q\rangle_{\oplus}$, where $\langle Q\rangle_{\oplus}$ is the full additive subcategory generated by $Q$. In fact, this is a decomposition of $\C$-bimodule categories, as can be seen from the next lemma.
       \begin{notation}\label{notn:decomp}
           For any $V\in \D$, we denote by $V = V_\C\oplus V_{\langle Q\rangle_{\oplus}}$ the corresponding summands $V_\C\in \C$ and $V_{\langle Q\rangle_{\oplus}}\in {\langle Q\rangle_{\oplus}}$.
       \end{notation}

    \comment{
    \begin{proof}
        \mbox{}
        \begin{enumerate}
            \item The tensor product of any two invertible objects in $\D$ is invertible. The objects in $\C$ are precisely those whose Jordan--H\"older composition factors are all invertible. Since the bifunctor $-\otimes-$ is biexact, it follows that the tensor product of any two objects in $\C$ has only invertible Jordan--H\"older composition factors, so it lies again in $\C$. This proves that $\C$ is a tensor subcategory of $\D$. 
            \item Let $M\in \D$, and set $m:=\dim \Hom_{\D}(Q, M)$. 
            
            Since $Q$ is injective and projective, there exists a unique (up to isomorphism) $M'\in \D$ such that $M\cong M'\oplus Q^{\oplus m}$ and $\dim \Hom_{\D}(Q, M')=0$. It remains to show that $M'\in \C$; that is, we need to show that the Jordan--H\"older composition factors of $M'$ are all invertible. Since $Q$ is injective and projective, it may appear as a composition factor in $M'$ if and only if it appears as a direct summand. Yet $\dim \Hom_{\D}(Q, M')=0$, implying that $Q$ does not appear as a composition factor in $M'$. This proves the statement.
        \end{enumerate}
    \end{proof}
}
The object $Q\otimes Q$ is projective (see \cite[Proposition 4.2.12]{EGNO}). Consider the decomposition of $Q\otimes Q$ into a direct sum of indecomposable projective objects.  
    For an isomorphism class of indecomposable projective objects $P\in \D$, denote by $[Q\otimes Q:P] $ the multiplicity of $P$ in the above decomposition.

    \begin{lemma}\label{lemma: tensor product of simple projective object}
        \mbox{}
        \begin{enumerate}
            \item Let $L \in \D$ be an invertible object, and let $P(L)$ denote its projective cover. Then $P(L)\in \C$. 
            \item $Q$ is self-dual, i.e. $Q^*\cong Q$.
            \item For any invertible object $L\in \D$, we have: $L\otimes Q\cong Q$.
            \item For any indecomposable projective $P\in \C$, we have $[Q\otimes Q:P] = 1$.
        \end{enumerate}
    \end{lemma}
    \begin{proof}
        \mbox{}
        \begin{enumerate}
            \item This follows directly from \cref{lem:blocks_almost_semsimple}(2).
            \item Since $Q$ is a non-invertible simple object in $\D$, so is $Q^*$. The assumption that $\D$ is a generalized near-group category now implies that $Q\cong Q^*$.
            \item Let $L\in \D$ be invertible, then $L\otimes -$ is an autoequivalence of $\D$ as an abelian category, with inverse $L^*\otimes-$, so $L\otimes Q$ is also simple and not invertible, thus $L\otimes Q\cong Q$.
            \item Let $L$ be the simple cosocle of $P$. We have:
            $$\dim\Hom_{\D}(Q,Q) =\dim\Hom_{\D}(Q,L\otimes Q) = \dim\Hom_{\D}(Q\otimes Q^*,L) = \dim\Hom_{\D}(Q\otimes Q ,L) = [Q\otimes Q:P].$$
            By Schur's lemma we have $$1=\dim \Hom_{\D}(Q,Q)=\dim \Hom_{\D}(Q\otimes Q ,L) =[Q\otimes Q:P].$$
        \end{enumerate}
    \end{proof}
    \comment{
    \begin{lemma}\label{cor:decomp_D_as_C_bimod}
        We have a decomposition of $\C$-bimodule categories $\D=\C\oplus \langle Q\rangle_{\oplus}$, where $\langle Q\rangle_{\oplus}$ is the full additive subcategory of $\D$ generated by $Q$.
    \end{lemma}
    }
    Let $r:=[Q\otimes Q:Q]$ denote the extension constant of $\D$. The above results imply that the Grothendieck ring $Gr(\D)$ depends only on $G:=\pic(\D)$, $r$ and $d:=\FPdim(P(\1))$. Note that $d$ equals the dimension of the projective cover of any invertible object, not just the unit.
    \begin{lemma}\label{cor:groth_ring}
        Let $\D$ be an almost-semisimple generalized near-group category, $G:=\pic(\D)$ its Picard group and $r$ its extension constant.
        Then the Grothendieck ring $Gr(\D)$ is a free $\Z$-module with $\Z$-basis $\{[L], [Q]~|~[L]\in G\}$. The multiplication in $Gr(\D)$ is given by $$[L][L']=[L \otimes L'], ~~[L][Q]=[Q][L]=[Q], ~~[Q]^2=r[Q]+d\sum_{[L]\in G} [L].$$ 
    \end{lemma}
    Similarly we obtain the following basic equation that we will use extensively:
    \begin{equation}\label{eq:q_square}
        q^2~=~rq+nd,
    \end{equation}
    where $q=\FPdim(Q)$, $n = |\pic(\D)|$, $d = \FPdim(P(\1))$ and $r = [Q\otimes Q:Q]$.

\section{Almost semisimple generalized near-group categories admitting a fiber functor}\label{sec: integral near-group categories}
    In this section, we consider almost semisimple generalized near-group categories admitting a fiber functor. For readability, from here on we will refer to such categories as \textit{admissible}.

The unitary fusion (semisimple) admissible categories have been classified by \cite{Iz, TY}. 

   From here on, we fix an admissible category $\D$, we denote by $\C\subset \D$ the (full) Serre subcategory generated by the invertible objects in $\D$, by $Q$ the non-invertible simple projective object in $\D$. 

    \subsection{Properties of admissible categories}
        Denoting the fiber functor $F: \D\to \Vect$, we set
        $H^*:=Coend(F)$, $ K^*:=Coend(F\rvert_{\C})$. Then $K^*\subset H^*$ is an embedding of Hopf algebras, and we have equivalences $$\D\cong H^*-Comod, \;\;\C\cong K^*-Comod.$$  Denote by $G$ the group of grouplike elements of $K^*$. Clearly, $K^*$ is a pointed Hopf algebra and $G \cong \pic(\C)$. We use the notation $H^*$ to avoid confusion with the next subsection, in which we consider the Hopf algebra $H = \operatorname{End}(F)$, which is dual to $H^*$. 

\begin{lemma}\label{lemma: decomposition as coalgebras}
        We have an isomorphism of coalgebras $H^*\cong K^*\oplus E(Q).$ Furthermore, $E(Q)$ has a structure of a $K^*$-bimodule, and the above isomorphism is an isomorphism of $K^*$-bimodules. 
    \end{lemma}
    \begin{proof}
        The isomorphism of coalgebras is obtained directly from the decomposition into $\C$-bimodule categories $\D=\C\oplus \langle Q\rangle_{\oplus}$. The fact that this is an isomorphism of bimodules follows from \cref{lem:E_X}.
    \end{proof}
    \begin{corollary}\label{crl: dim(K) divides q^2}
        We have: $\dim K^* ~\rvert~ \FPdim(Q)^2$.
    \end{corollary}
    \begin{proof}
        The algebra $H^*$ is free over $K^*$ by the Nichols-Zoeller Theorem (\cite{NZ}), so $\dim K^* ~\rvert~\dim H^*$. By \cref{lemma: decomposition as coalgebras}, $\dim(H^*) = \dim(K^*)+\dim(E(Q))$, thus $\dim K^* ~\rvert~\dim E(Q)=\FPdim(Q)^2$.
    \end{proof}
Together with \cref{eq:q_square}, \cref{crl: dim(K) divides q^2} implies:
    \begin{corollary}\label{cor:divisibility}
        Let $q:=\FPdim(Q)$, $d:=\FPdim (P(\1))$ and $n:=|G|$. Then $q|(nd)$ and $(nd)|q^2$. Furthermore, $r$ is divisible by $\frac{nd}{q}$.
    \end{corollary}
    \begin{lemma}\label{lem:projective_G_rep}
        There is a natural structure of a projective action $G\to PGL(F(Q))$ which is compatible with the $G$-action on $E(Q)\cong F(Q)\otimes F(Q)^*\subset H^*$ by left multiplication.
    \end{lemma}
    \begin{proof}
        The $G$-action on $E(Q)\cong F(Q)\otimes F(Q)^*\subset H^*$ by left multiplication makes $E(Q)$ into a coalgebra in the category of $G$-modules. Dualizing, we get that $\End_{\kk}(F(Q))$ is an algebra in the category of $G$-modules. Since $\End_{\kk}(F(Q))$ is a simple algebra, we have (e.g. by Skolem-Noether theorem) that the action of $G$ on $\End_{\kk}(F(Q))$ is inner, inducing a compatible projective action $G\to PGL(F(Q))$.
    \end{proof}

    \begin{prop}\label{prop: dim(Q) in non-semisimple integral near-group categories}
        Assume that $H^2(G,\kk^\times) = 0$ (e.g. $G$ is cyclic). Then $|G|$ divides $\FPdim(Q)$.
    \end{prop}
    \begin{proof}
    Since $H^2(G,\kk^\times) = 0$, the projective-linear action of $G$ on $F(Q)$ lifts to a linear action of $G$ on $F(Q)$.
    
    Consider the embeddings of algebras $\kk G \subset K^*\subset H^*$. By the Nichols-Zoeller theorem, the algebras $K^*$ and $H^*$ are free over the Hopf subalgebra $\kk G$. So when considered as left $G$-modules, both $K^*, H^*$ are isomorphic to direct sums of several copies of the left regular representation $\kk G\cong \bigoplus_{\pi \in Irr(G)} \pi^{\oplus \dim \pi}$. 
    
    Recall the decomposition of $K^*$-modules $H^*=K^*\oplus E(Q)$. Restricting the $K^*$-action to that of $G$, we conclude that $E(Q)$ is a direct sum of several copies of the left regular $G$-representation $\kk G$, so there exists $m> 0$ such that $E(Q)\cong \left(\bigoplus_{\pi \in Irr(G)} \pi^{\oplus \dim \pi}\right)^{\oplus m}$ as $G$-modules.

    On the other hand, we have an isomorphism of $G$-modules $E(Q)\cong F(Q)\otimes F(Q)^*$ due to the compatibility stated in \cref{lem:projective_G_rep}. 

  Let $F(Q)\cong  \bigoplus_{\pi \in Irr(G)} \pi^{\oplus m_{\pi}}$ be the decomposition of $F(Q)$ into irreducible $G$-modules. 

  Then the multiplicity of the trivial $G$-module $\1$ in $F(Q)\otimes F(Q)^*$ is $\sum_{
  \pi \in Irr(G)} m_{\pi}^2$, while we saw that the multiplicity of $\1$ in the $G$-module $E(Q)$ is $m$. So $m=\sum_{
  \pi \in Irr(G)} m_{\pi}^2$. Now the dimension of $E(Q)$ is $$ m\sum_{  \pi \in Irr(G)} \dim(\pi)^2=  \left(\sum_{  \pi \in Irr(G)} m_\pi^2 \right)\left(\sum_{  \pi \in Irr(G)} \dim(\pi)^2 \right)=(\dim F(Q))^2=\left(\sum_{ \pi \in Irr(G)} m_\pi\dim(\pi)\right)^2$$

  By the Cauchy--Schwarz inequality, we conclude that the sequence $(m_{\pi})_{\pi \in Irr(G)}$ is a multiple of the sequence $(\dim{\pi})_{\pi \in Irr(G)}$, with the ratio $m':=m_{\1}$. Thus as a left $G$-module, $F(Q) \cong \left(\kk G\right)^{\oplus m'}$, and $\dim F(Q)=m'|G|$, as required.
\end{proof}
\begin{crl}\label{corollary: q-r divides d in trivial central extensions}
   Let $q:=\FPdim(Q)$, $d:=\FPdim (P(\1))$ and $n:=|G|$. Assume $H^2(G,\kk^\times) = 0$. 
     
     Then $n|q$ and $d =\frac{q}{n}(q-r)$. In particular, $q-r$ and $\frac{q}{n}$ divide $d$.
\end{crl}
\begin{proof}
    By \cref{prop: dim(Q) in non-semisimple integral near-group categories} we know that $n\rvert q$. Dividing \cref{eq:q_square} by $n$ gives the required result.
\end{proof}

    \subsection{The structure of the Hopf algebra of flip-compatible admissible category}
        In this section, we fix an admissible category $\D$, with a fiber functor $F:\D\to \Vect $. Let $H := \End(F)$ and $K:=\End(F|_\C)$ be the corresponding Hopf algebras. We denote by $\pi:H\to K$ the corresponding Hopf algebra homomorphism.

        Fix an equivalence $\D\cong\Rep(H)$ whose composition with the forgetful functor $\Rep(H)\to \Vect$ gives $F$. We will identify $\D$ with $\Rep(H)$ under this equivalence.

        The decomposition $\D=\C\oplus {\langle Q\rangle_{\oplus}}$ gives a pair of orthogonal central idempotents $e,(1-e)$ in $H$, so that for any $V\in \Rep(H)$, we have $$e(V)=V_{\C},\;\; (1-e)(V)=V_{\langle Q\rangle_{\oplus}}$$ (see \cref{notn:decomp}). 
        We have an isomorphism of algebras $He\cong K$, while $H(1-e)$ is a semisimple algebra. 
        
        Denote by $\tau_{V,W}$ the flip linear map $\tau_{V,W}:V\otimes W\to W\otimes V$ for $V,W\in \Rep(H)$. In general, this is not a map of $H$-modules, but for $h\in H$ we have an equality in $\Hom_{\kk}(V\otimes W, W\otimes V)$:
        \begin{equation}\label{eq: flip-compatability}
            \tau_{V,W}\circ \Delta(h)  = \Delta^{op}(h)\circ \tau_{V,W}
        \end{equation}
        where $ \Delta^{op}:=\tau_{H, H}\circ \Delta$.

\begin{lemma}\label{lem:properties_of_e}
    The idempotent $e$ satisfies the following properties: 
    \begin{enumerate}
              \item $\varepsilon(e)=1$ and $S(e)=e$.  
              \item $e\in H^{co \pi} \cap \prescr{co\pi}{H}$, and $eh=he=\varepsilon(h)e$ for any $h\in  H^{co \pi} \cap \prescr{co\pi}{H}$.
\item $\ker\pi=H(1-e)$.
    \end{enumerate}
\end{lemma}
\begin{proof}
\mbox{}
\begin{enumerate}
\item The trivial module $\1_\D$ lies in $\mathcal{C}$, so $e$ acts on it as the identity, implying $\varepsilon(e)=1$. The subcategory $\mathcal{C}$ is closed under taking duals, as is $\langle Q\rangle_{\oplus}$, so $S(e)=e$.
        \item
Recall that $\langle Q\rangle_{\oplus}\otimes\mathcal{C}\subseteq\langle Q\rangle_{\oplus}$ and $\mathcal{C}\otimes\langle Q\rangle_{\oplus}\subseteq\langle Q\rangle_{\oplus}$.
Therefore, for any $V\in\mathcal{D}$ and $U\in \mathcal{C}$, $$(V\otimes U)_{\mathcal{C}}=V_{\mathcal{C}}\otimes U.$$
So $\Delta(e) \rvert_{V\otimes U}=(e\otimes1) \rvert_{V\otimes U}$. 
Take $V:=H$ and $U:=K$, with $H$ acting on $\mathit{K}$ through $\pi$. Both elements $\Delta(e),~ e\otimes 1$ act by left multiplication on $H\otimes K$. Applying these to $1_H\otimes1_K \in H\otimes K$ gives $(\mathrm{id}\otimes\pi)\Delta(e)=e\otimes1$, so $e\in H^{co \pi}$. 
The mirror argument with $U\otimes V$ gives $(\pi\otimes\mathrm{id})\Delta(e)=1\otimes e$, so $e\in  \prescr{co\pi}{H}$.

Let $h\in H^{co \pi} \cap \prescr{co\pi}{H}$. Apply $\pi\otimes\mathrm{id}$ to $(\mathrm{id}\otimes\pi)\Delta(h)=h\otimes1$. This gives $\Delta_K(\pi(h))=\pi(h)\otimes1$.
Applying $\varepsilon\otimes\mathrm{id}$ then gives $\pi(h)=\varepsilon(h)1$.
Hence $\pi(eh)=\varepsilon(h)\pi(e)=\pi(\varepsilon(h)e)$, since $\pi(e)=\varepsilon(e)1=1$ by (1).
The restriction of $\pi$ to $eH=He$ is an isomorphism onto $\mathit{K}$. So $eh=\varepsilon(h)e$. Since $\mathit{e}$ is central, also $he=\varepsilon(h)e$.

       \item This follows from the decomposition $H=He\times H(1-e) \cong K\times H(1-e)$.
       
    \end{enumerate}
\end{proof}
\begin{prop}\label{prop: exact sequence of Hopf algebras}
    The following conditions are equivalent: 
    \begin{enumerate}
        \item For all $V,W\in \D$, the flip morphism $\tau_{V, W}$ sends the subspace $F((V\otimes W)_\C) \subset F(V)\otimes F(W)$ into the subspace $F((W\otimes V)_{\C}) \subset F(W)\otimes F(V)$.
        \item For all $V, W\in \Rep(H)$, we have: $$\tau_{W, V}\left((W\otimes V)_{\mathcal{C}}\right)=(V\otimes W)_{\C} \;\;\text{ and } \;\;(V\otimes W)_{\langle Q\rangle_{\oplus}}=\tau_{W, V}\left((W\otimes V)_{\langle Q\rangle_{\oplus}}\right)$$ 
        \item $\Delta(e)=\Delta(e)^{op}$, i.e. $\Delta(e)\in H\otimes H$ is a symmetric tensor.
    
        \item The homomorphism $\pi: H\to K$ is normal.
    \end{enumerate}   

If these conditions hold, then $H^{co \pi}<H$ is a normal semisimple Hopf subalgebra and $e\in H^{co \pi}$.
\end{prop}

          In particular, any of the conditions (1)-(3) in the above proposition implies that there exists an exact sequence of Hopf algebras
        $$
            \kk\to H^{co \pi}\to H\xrightarrow{\pi} K\to \kk,
        $$
        and $H$ is a crossed product of $H^{co \pi}$ with $K$.

    \begin{proof}\mbox{}
    
    \begin{enumerate}
    \item[\textbf{[(1) $\Longleftrightarrow$ (2)]}] Condition (2) clearly implies that $\tau_{W, V}\left((W\otimes V)_{\mathcal{C}}\right)\subset (V\otimes W)_{\C}$, which is precisely condition (1). In fact, the implication in the other direction is true as well. Assume that (1) holds. The space $(V\otimes W)_{\langle Q\rangle_{\oplus}}$ is the annihilator of $({W^*}\otimes {V^*})_\C$ under the evaluation pairing, because $\C$ and ${\langle Q\rangle_{\oplus}}$ are each closed under duality. The flips preserve the pairing, so applying (1) to the pair $({W^*},{V^*})$ yields (2).
        \item[\textbf{[(2) $\Longleftrightarrow$ (3)]}] Fix $V, W\in \Rep(H)$. By definition, $\Delta(e) $ acts on $V\otimes W$ by projection onto the direct summand $(V\otimes W)_{\C}$ with kernel $(V\otimes W)_{\langle Q\rangle_{\oplus}}$. Meanwhile, $\Delta^{op}(e)$ is the idempotent on $V\otimes W$ with image $\tau_{W, V}\left((W\otimes V)_{\mathcal{C}}\right)$ and kernel $\tau_{W, V}\left((W\otimes V)_{\langle Q\rangle_{\oplus}}\right)$. 
        Applying these to the left $H$-module $H$, we get: $ \Delta(e)=\Delta(e)^{op}$ if and only if $\tau_{H, H}\left((H\otimes H)_{\C}\right)=(H\otimes H)_{\C}$ and $\tau_{H, H}\left((H\otimes H)_{\langle Q\rangle_{\oplus}}\right)=(H\otimes H)_{\langle Q\rangle_{\oplus}}$.
        
        So $ \Delta(e)=\Delta(e)^{op}$ if and only if for any $V, W\in \Rep(H)$, we have: $\tau_{W, V}\left((W\otimes V)_{\mathcal{C}}\right)=(V\otimes W)_{\C}$ and $(V\otimes W)_{\langle Q\rangle_{\oplus}}=\tau_{W, V}\left((W\otimes V)_{\langle Q\rangle_{\oplus}}\right)$.

        \item[\textbf{[(3) $\Longrightarrow$ (4)]}] Assume that $\Delta(e)=\Delta(e)^{op}$. 
        We will use \cref{lem:normal_equiv_cond} to establish the existence of a suitable subalgebra $A$ as in \cref{lem:normal_equiv_cond}.
        
        Let us write $\Delta(e) =: \sum_i a_i\otimes b_i$ so that $\{a_i\}_i$ and $\{b_i\}_i$ are two sets of linearly independent vectors in $H$, and let $R$ be the vector space spanned by $\{a_i\}_i$ (note that $R$ does not depend on the choice of presentation $\Delta(e) = \sum_i a_i\otimes b_i$). By our assumption, $\Delta(e)=\Delta(e)^{op}= \sum_i b_i\otimes a_i$ so 
        $$R=\mathrm{span}\{a_i\}_i=\mathrm{span}\{b_i\}_i.$$
        Let us consider the subalgebra $A$ of $H$ generated by $R$.
        
Coassociativity of $\Delta$ implies that $\Delta(R)\subset H\otimes R\cap R\otimes H= R\otimes R$, so $A$ is in fact a bialgebra. Thus $A$ is a sub-bialgebra of the finite-dimensional Hopf algebra $H$, hence a Hopf subalgebra.

Recall that by \cref{lem:properties_of_e}, we have: $e\in H^{co \pi}\cap \prescr{co \pi}{H}$. Since $H^{co \pi}$ is a left coideal and $\prescr{co \pi}{H}$ is a right coideal, we get that $b_i\in H^{co \pi}$, $a_i\in \prescr{co \pi}{H}$ for every $i$, so $R\subset H^{co \pi}\cap \prescr{co \pi}{H}$. Since $H^{co \pi}, \prescr{co \pi}{H}$ are subalgebras, we conclude: $A\subset H^{co \pi}\cap \prescr{co \pi}{H}$. By \cref{lem:coinvar_kernel}, we have: $\pi(h) = \varepsilon(h)1_K$ for any $h\in A$.


Next, $e\in A$ since
$$e=(\id\otimes \varepsilon)\Delta(e)=\sum_i a_i\varepsilon(b_i) \in \mathrm{span}\{a_i\}_i =R\subset A.$$
So $(1-e)\in A$. By \cref{lem:properties_of_e}(3), $H(1-e) = \ker \pi$ so we conclude that $\ker \pi \subset HA^+ \subset \ker \pi$ and thus  $\ker \pi = HA^+$.
\comment{
By the Nichols-Zoeller Theorem (\cite{NZ}), $H$ is free over $A$. Let $m:=  rk_A H=\frac{\dim(H)}{\dim(A)} \in \Z$; we will show that $m=\dim K$.

Indeed, consider the short exact sequence of left $A$-modules
        $$
            0\to A^+\to A\xrightarrow{\varepsilon} \1_A\to 0.
        $$ where $\1_A $ stands for the trivial $1$-dimensional $A$-module.
        Applying $H\otimes_A-$ gives a short exact sequence of $H$-modules
        $$
            0\to H A^+\to H\to H\otimes_A \1_A\to 0
        $$
        and $\dim H\otimes_A \1_A = rk_A H=m$.
     We conclude that $$ m = \frac{\dim(H)}{\dim(A)}=\dim H-\dim HA^+ = \dim H-\dim \ker \pi = \dim K$$ as required. 
     }
     Hence $\pi$ is normal.
     \item[\textbf{[(4) $\Longrightarrow$ (3)]}]  Assume that $\pi$ is normal, so $ A:=H^{co \pi}=\prescr{co \pi}{H}$ is a Hopf subalgebra. By \cref{lem:properties_of_e}, we have: $e\in A$, and in particular, $\Delta(e), \Delta^{op}(e) \in A\otimes A$. Thus we need to check that $\Delta_A(e)= \Delta^{op}_A(e)$ where $\Delta_A$ is comultiplication in $A$.

    Recall from \cref{lem:properties_of_e} that $e$ is both a left and a right integral in $H^{co \pi}$. Since $\eps(e)\ne 0$, we conclude (see \cite{LS} or \cite[Proposition 6.5.6.]{EGNO}) that $H^{co \pi}$ is semisimple.  
    Since $A$ is a finite-dimensional semisimple Hopf algebra over a field of characteristic zero, it is also cosemisimple.
By \cite[Proposition 7.18.14]{EGNO} (due to Larson and Radford), if $A$ is a semisimple and cosemisimple Hopf algebra then $\Delta_A(I)= \Delta^{op}_A(I)$ for any left integral $I\in A$. Thus $\Delta_A(e)= \Delta^{op}_A(e)$, as required.

        \end{enumerate}
Thus we proved that conditions (1)-(4) are equivalent and imply that $A$ is semisimple. By \cref{lem:normal_equiv_cond}, $A=H^{co \pi}=\prescr{co \pi}{H}$ is a normal Hopf subalgebra of $H$.

       \end{proof}

    \begin{definition}\label{def: flip-compatible fiber functor}
        We say that a fiber functor $F:\D\to \Vect$ is \textit{flip compatible} if one of the conditions in \cref{prop: exact sequence of Hopf algebras} holds. 
        
        Given a pair $(\D, F:\D\to \Vect)$ consisting of an admissible category $\D$ with a flip-compatible fiber functor, we write that this pair is \textit{FCA} for short.
    \end{definition}  
    
    \begin{rem}\label{rmk:flip_compatible}
    
    \mbox{}
    
    \begin{enumerate}
        \item The fiber functor $F:\D\to \Vect$ is flip-compatible if and only if the condition in \cref{def: flip-compatible fiber functor} holds for $V=W=Q$. In particular, for $r=0$ every fiber functor is flip-compatible.

        \item By \cite{EG2}, the only family of unitary near-group categories admitting a fiber functor with $r>0$ are $\Rep(\mathbb F_q\rtimes \mathbb F_q^\times)$ which are Tannakian, so the forgetful functor in this case is flip-compatible. We emphasize that this relies on the classification of unitary near-group fusion categories. We are not aware of a general categorical reason for compatibility in the semisimple case.
    \end{enumerate}
       
    \end{rem}
        \begin{conj}
        For every almost semisimple generalized near-group category admitting a fiber functor $\D$, there exists a fiber functor $F:\D\to \Vect$ which is flip-compatible.
    \end{conj} 

        From now on, fix an FCA pair $(\D, F:\D\to \Vect)$ and let        $$
            \kk\to H^{co \pi}\to H\xrightarrow{\pi} K\to \kk,
        $$ be the short exact sequence arising from \cref{prop: exact sequence of Hopf algebras}. 
        We will now list some further properties of the Hopf subalgebra $A:=H^{co \pi}$.

    Let $Irr(A)$ denote the set of isomorphism classes of simple $A$-modules.
    \begin{prop}\label{prop: structure of Res(Q)}

    We have an isomorphism of left $A$-modules
        $$Q \cong  \left(\bigoplus_{L\in Irr(A)\setminus\{\1_A\}}L^{\oplus \dim(L)}\right)^{\oplus (q-r)}.$$           
    \end{prop}
  
    \begin{proof}
    Consider the left $H$-module $H$ and its quotient, the left $H$-module $K$. Since $\pi(h)=\varepsilon(h)1$ for $h\in A$, the module $K$ is a trivial $A$-module of dimension $ \dim(K)$.

    
        By Frobenius reciprocity we have: for any $L\in Irr(A)$, $L\neq \1_A$ we have:
        $$
           \dim\Hom_H(H\otimes_A L, K)=\dim\Hom_A(L, K) =0
        $$
        so $H\otimes_A L = Q^{\oplus {\alpha_L}}$ for some $\alpha_L \in \Z_{\geq 0}$.      
        Since $A$ is semisimple, we have a decomposition of the $A$-module $Q$ into a direct sum of irreducible $A$-modules. By Frobenius reciprocity, this decomposition is given by $$Q\cong \bigoplus_{L\in Irr(A)\setminus\{\1_A\}}L^{\oplus \alpha_L}.$$
        We need to prove that $\alpha_L=\dim(L)(q-r)$ for any $L\in Irr(A)\setminus\{\1_A\}$.
        By the Nichols-Zoeller theorem (\cite{NZ}), $H$ is a free left $A$-module of rank $\dim(K) = nd$. So $\dim H\otimes_A L   = nd\cdot \dim(L)$. On the other hand, $$\dim(H\otimes_A L) = \alpha_L \dim Q = \alpha_L q,$$ so $\alpha_L = \frac{nd}{q}\dim(L)$. By \cref{eq:q_square} we get $q-r = \frac{nd}{q}$, completing the proof.
    \end{proof}
    \begin{remark}
        We already know by \cref{crl: dim(K) divides q^2} and \cref{eq:q_square} that $(q-r )\rvert q$. \cref{prop: structure of Res(Q)} gives a module-theoretic interpretation of this divisibility.
    \end{remark}

    \begin{prop}\label{prop: A is a minimal normal Hopf subalgebra}
        
        For every non-trivial normal Hopf subalgebra $N\subset H$, we have: $N$ acts non-trivially on $Q$, and $A\subset N$.
    
    \end{prop}
    \begin{proof}
     Let $N\subset H$ be a normal Hopf subalgebra.
       Let $\C'$ be the full subcategory of $\Rep(H)$ whose objects are $H$-modules on which the restricted action of $N$ is trivial. This subcategory is a tensor subcategory (a subcategory closed under taking subquotients, direct sums and tensor products).

        Consider the $H$-module $B:=H/HN^+$. Clearly, $B\in \C'$. 
        
        If $Q\in \C'$ then $Q\otimes Q\in \C'$. Yet any indecomposable projective $H$-module which is not isomorphic to $Q$ occurs as a direct summand in $Q\otimes Q$ (see \cref{lemma: tensor product of simple projective object}), so in this case $\C'$ contains all the projective $H$-modules and so $\C'=\Rep(H)$. This would imply that $N = \kk$. 
        
        Now assume that $Q\notin \C'$. Then $\C'\subseteq \C$ so the idempotent $(1-e)$ acts trivially on any $M\in \C'$, and in particular on $B$. Since $\ker \pi = H(1-e)$, we conclude that the quotient map $\pi':H\twoheadrightarrow B$ factors through $\pi$. Since $N$ is normal, this implies that $A=H^{co \pi}\subseteq H^{co \pi'} = N$ (the last equality follows from \cref{lem:normal_equiv_cond}).
    \end{proof}

    Using the results above, we can prove our main theorem (see \ref{main theorem}):
    
    \begin{theorem}\label{thrm:main_theorem_in_full}
        Let $(\D,F:\D\to \Vect)$ be an FCA pair, and let $\C\subset \D$ be the (full) Serre subcategory generated by the invertible objects. Let $H:=\End(F)$ and $K := \End(F|_{\C})$. We have an algebra isomorphism $$H \cong A\#_\sigma K$$ where $A\subset H$ denotes the smallest non-trivial normal Hopf subalgebra of $H$ and $\sigma:K\otimes K\to A$ is a $2$-cocycle. Moreover, $A$ is semisimple.

        Equivalently, there exists a canonical exact sequence of tensor categories:
        $$
            \C\to \D\to \E
        $$
        where
        $\E:=\Rep(A)$ is a fusion category.

    \end{theorem}

    \begin{proof}
        By \cref{prop: exact sequence of Hopf algebras} we have an exact sequence 
        $$
            \kk \to A\to H\to K\to \kk.
        $$
        for a semisimple Hopf subalgebra $A\subset H$. By \cref{prop: A is a minimal normal Hopf subalgebra}, $A$ is the minimal non-trivial normal Hopf subalgebra of $H$.
        
        By \cite{Sch} (see also \cref{ssec:crossed_prod}), $H$ is a crossed product of $A$ and of $K$ by some $2$-cocycle $\sigma$. The existence of the exact sequence of tensor categories follows from \cite{BN} (see also \cref{ssec:exact_seq_tensor_cats}), and the category $\E$ is fusion because $A$ is semisimple.
    \end{proof}
    \begin{prop}\label{prop: A is central when dim(A)=2}
        If $\dim(A) = 2$ then $A\subset Z(H)$.
    \end{prop}
    \begin{proof}
        By \cref{prop: exact sequence of Hopf algebras}, the central idempotent $e\in H$ defined in the beginning of this section lies in $A$. By our assumption, $\dim A=2$ so $A = \mathrm{span}\{1,e\}\subset Z(H)$ as required. 
    \end{proof}
    \begin{thm}\label{thm: integral almost semisimple with r=0}
        Let $(\D,F:\D\to \Vect)$ be an FCA pair, and let $\C\subset \D$ be the (full) Serre subcategory generated by the invertible objects. Let $H:=\End(F)$ and $K := \End(F|_{\C})$. 
        
        The following conditions are equivalent:
        \begin{enumerate}
            \item $\D$ has extension constant $r:=[Q\otimes Q:Q]=0$.
            \item We have an isomorphism $H\cong \kk[\Z/2\Z]\#_\sigma K$ for some $2$-cocycle $\sigma:K\otimes K\to \kk[\Z/2\Z]$ and the image of $\kk[\Z/2\Z]$ in $H$ is central.
        \end{enumerate}
    \end{thm}
    \begin{proof}
        Let $Q$ be the unique non-invertible simple object in $\D$, with $q:=\dim(Q)$ as before. 
    
        \begin{enumerate}
            \item[\textbf{[(2) $\Longrightarrow$ (1)]}] Assume that $H\cong \kk[\Z/2\Z]\#_\sigma K$. Comparing the dimensions, we get: $\dim(H) = 2\dim(K)$, so $$nd+q^2 = \dim(H) = 2\dim(K) = 2nd$$ which implies that $q^2=nd$. \cref{eq:q_square} now implies that $r=0$.
            \item[\textbf{[(1) $\Longrightarrow$ (2)]}] Assume that $\D$ has generalized fusion rule $(G,0)$, i.e. $r=0$.
            
  By \cref{thrm:main_theorem_in_full}, we have a normal semisimple Hopf subalgebra $A\subset H$ such that $H\cong A\#_\sigma K$ for some $2$-cocycle $\sigma:K\otimes K\to A$. 
  
            By \cref{prop: structure of Res(Q)}, we have an isomorphism of left $A$-modules
        $$Q \cong  \left(\bigoplus_{L\in Irr(A)\setminus\{\1_A\}}L^{\oplus \dim(L)}\right)^{\oplus q}.$$
            Since $\dim Q=q$, we get that $\dim \bigoplus_{L\in Irr(A)\setminus\{\1_A\}}L^{\oplus \dim(L)} =1$, so $|Irr(A)|=2$ and $\dim(L)=1$ for any $L\in Irr(A)$. Thus $\dim A=2$ and since $A$ is semisimple, we conclude that $A\cong \kk[\Z/2\Z]$. By \cref{prop: A is central when dim(A)=2}, $A\subset Z(H)$ which completes the proof.
          
        \end{enumerate}

        
    \end{proof}

\section{Extensions of Taft algebras}\label{sec: Extending H_{n,d}}
    \subsection{Taft algebras}\label{ssec:gen_Taft_alg}
        Let $ n\in \N_{>1}$.
        Let $\omega \in \kk$ be a primitive $n$-th root of unity. The Taft algebra is a Hopf algebra, whose underlying algebra is defined as $$H_{n}=\langle\, g,x\,|\,g^n=1, x^n=0, xg=\omega gx\,\rangle.$$
        The coalgebra structure $(\Delta, \epsilon)$ and an antipode $S$ are given by
        \begin{align*}
            & \Delta(g) = g\otimes g,~\epsilon(g) = 1\\
            & \Delta(x) = g \otimes x + x\otimes 1,~\epsilon(x) = 0 \\
            & S(g) = g^{-1},~S(x) = -g^{-1}x.
        \end{align*}

        \begin{example}
            The algebra $H_{2}$ is Sweedler's Hopf algebra.
        \end{example}

        It is well-known that the category $H_n-mod$ is a pointed non-semisimple finite tensor category with $n$ simple $1$-dimensional objects (up to isomorphism), each of which has a projective cover of dimension $n$ (see for instance, \cite{T1, EO}). We are particularly interested in the Taft algebra due to the following property:
        \begin{lemma}
            The Picard group of $H_n-mod$ is the cyclic group of order $n$.
        \end{lemma}

        The primitive idempotents of $H_{n}$ are given by $e_i = \frac{1}{n}\sum_{j=0}^{n-1}\omega^{-ji}g^j$ for $i\in \Z/n\Z$. The dual to $L_i$ is $L_i^*\cong L_{n-i}$. The simple object $L_i$ corresponding to $e_i$ is the $1$-dimensional module generated by an element $u_i$ with the actions $g. u_i = \omega^i u_i,~x. u_i = 0$. The indecomposable projective objects in $H_n-mod$ are $P_i =H_{n} e_i$ (projective cover of $L_i$) for $i\in \Z/n\Z$. The module $P_i$ has a natural basis $ \{e_i,...,x^{n-1}e_i\}$ and $$socP_i=span\{x^{n-1}e_i\} \cong L_{i-(n-1)}.$$ 
        
        The socle and radical filtrations of the projective objects coincide, and: 
        $$\frac{rad^j P_i}{rad^{j+1} P_i} \cong L_{i-j} \;\; \text{ for } \;\; 0\leq j< n.$$

    \subsection{Extensions of Taft algebras}

        In this section we will define a new Hopf algebra $A_n$ which will have the Taft algebra $H_{n}$ as its quotient.
        Let $$A_{n}:=\mathbbm{k}\langle \,a,z\,|\,a^{2n}=1, z^n=a^n-1, za=\omega az\,\rangle.$$
   
        This is a finite dimensional algebra; we will now describe a Hopf algebra structure on $A_{n}$. 
    
        We define a coalgebra structure on $A_{n}$ by defining the comultiplication on the generators in the following way:
        $$
            \Delta(a)=a\otimes a,~\epsilon(a)=1, \\
            ~\Delta(z) = a\otimes z +z\otimes 1,~ \epsilon(z)=0
        $$
        We then extend the comultiplication and counit morphisms to the entire algebra $A_{n}$ by making $\Delta$ and $\epsilon$ algebra morphisms. Define the antipode map $S$ to be $S(a)=a^{-1}, S(z) = -a^{-1}z$. 

        \begin{lemma}\label{lemma: delta z^r}
            For every $k\geq 0$, we have:
            $$
                \Delta(z^k) = \sum_{i=0}^{k} \theta^{(k)}_i a^{i}z^{k-i}\otimes z^i
            $$
            where $\theta^{(k)}_i \in \mathbbm{k}$ and $ \theta^{(k)}_0= \theta^{(k)}_k=1$. Furthermore, $\theta^{(n)}_i=0$ for $1\leq i\leq n-1$, so $$\Delta(z^n)=z^n\otimes 1+a^{n}\otimes z^n.$$
        \end{lemma}
        \begin{proof}
            We shall prove by induction on $k$ that
            $$  \Delta(z^k) = \sum_{i=0}^{k} \theta^{(k)}_i a^{i}z^{k-i}\otimes z^i$$ for some constants $\theta^{(k)}_i\in \mathbbm{k}$, and derive a recurrent formula for these constants. Indeed, for $k=0$ we have $\Delta(z^0)=1$ so $\theta^{(0)}_0=1$. Assuming the claim holds for some $k\geq 0$, we have:
            \begin{align*}     
                \Delta(z^{k+1}) &= (a\otimes z+z\otimes 1)\left(\sum_{i=0}^{k} \theta^{(k)}_i a^{i}z^{k-i}\otimes z^i\right) = \sum_{i=0}^{k} \theta^{(k)}_i a^{(i+1)}z^{k-i}\otimes z^{i+1} +\\&+\sum_{i=0}^{k} \theta^{(k)}_i za^{i}z^{k-i}\otimes z^{i} = \sum_{i=0}^{k} \theta^{(k)}_i a^{(i+1)}z^{k-i}\otimes z^{i+1} +\sum_{i=0}^{k} \theta^{(k)}_i \omega^{i} a^{i}z^{k+1-i}\otimes z^{i} \\
                &=  \sum_{i=0}^{k+1} \left(\theta^{(k)}_{i-1}+\theta^{(k)}_i \omega^{i}\right) a^{i}z^{k+1-i}\otimes z^{i}  .
            \end{align*}
            where we denote $\theta^{(k)}_{-1}:=0$. Let us also denote $\theta^{(k)}_i:=0$ for $i>k$.

            Thus for each $k\geq 1$ we have: $\theta^{(k)}_0=1$ and  $\theta^{(k)}_i= \theta^{(k-1)}_{i-1}+\theta^{(k-1)}_i \omega^{i}$ for $i\geq 1$. Taking a formal parameter $q$, consider the polynomials $\theta^{(k)}_i\in \Z[q]$ given by $\theta^{(k)}_0(q)=1$ and  $\theta^{(k)}_i(q)= \theta^{(k-1)}_{i-1}(q)+\theta^{(k-1)}_i(q) q^{i}$. The latter is known as the $q$-Pascal identity, and
           the polynomial $\theta^{(k)}_i$ is the $q$-binomial coefficient $$\binom{k}{i}_q=\frac{[k]_q!}{[i]_q! [k-i]_q!} \text{ where } [m]_q:=\frac{q^m-1}{q-1}.$$
            Specializing at $q:=\omega$ gives:   $\theta^{(k)}_0=\theta^{(k)}_k=1$ for all $k\geq 0$ and $\theta^{(n)}_i=0$ for all $1\leq i\leq n-1$. This completes the proof of the lemma.
        \end{proof}
        \begin{lemma}
            The maps $\Delta, \epsilon, S$ define a Hopf algebra structure on $A_{n}$.
        \end{lemma}
        \begin{proof}
            Let us show that the maps $\Delta, \epsilon, S$ are well defined; that is, we need to show that $\Delta,~\epsilon, ~S$ respect the defining relations in the algebra $A_{n}$.  

            A direct calculation shows that the relation $za=\omega az$ is preserved by both $\Delta$ and $\epsilon$. The calculation below shows that this relation is respected by $S$ as well:
            \begin{align*}
                &S(az)=S(z)S(a)=-a^{-1}za^{-1} = -a^{-1}za^{2n-1} = -(\omega^{-1})a^{-2}z\\
                &S(za) = S(a)S(z)=a^{-1}(-a^{-1}z) = -a^{-2}z,
            \end{align*}
            and since $S$ is $\kk$-linear we obtain that $S(az)=S(\omega^{-1} za)$. Next, by \cref{lemma: delta z^r} we have:
            $$
                \Delta(z^n) = a^{n}\otimes z^n+z^n\otimes 1 = a^n\otimes (a^n-1)+(a^n-1)\otimes 1 = a^n\otimes a^n-1\otimes 1 = \Delta(a^n-1),
            $$
            which implies that $\Delta$ respects the relation $z^n=a^n-1$. It is easy to see that $\epsilon$ respects this relation as well. Using induction on $k$, we immediately obtain:
            $$
                S(z^k)=S(z)^k = (-1)^k\omega^{\sum_{i=1}^{k-1}-i}a^{-k}z^k=(-1)^k\omega^{-\frac{k(k-1)}{2}}a^{-k}z^k.
            $$
           Denoting by $\mu$ a primitive $2n$-th root of unity such that $\mu^2 = \omega$, we have for $k=n$:
            \begin{align*}
                S(z^n) &= (-1)^na^{-n}z^n\mu^{-n(n-1)}= a^{-n}z^n(-1)^n\mu^{-n(n-1)}=-a^{-n}z^n = (-1)(1-a^{-n}) = a^{-n}-1\\
                &=S(a^n-1).
            \end{align*}
            The second equality follows from the fact that $\mu^n=-1$.
        
            A similarly simple calculation shows that the maps $\Delta,\epsilon,S$ satisfy the conditions making $A_{n}$ into a Hopf algebra.
            \end{proof}

        The following proposition is straightforward:
        \begin{prop}\label{prop:A_n_n_surjects_onto_H_n_n}
         The map $\psi:A_{n}\to{H_{n}}$ given by $\psi(a):=g,\psi(z):=x$ is a surjective Hopf algebra morphism.
        \end{prop}

        In particular, this implies that $A_n-mod$ contains $H_n-mod$ as a full tensor subcategory.

        \begin{prop}\label{prop: A_n is a crossed product}
            The following is an exact sequence of Hopf algebras:
            $$
                \kk \to \kk[\Z/2\Z]\to A_n\to H_n\to \kk
            $$
            in other words, there exists a $2$-cocycle $\sigma:H_n\otimes H_n \to \kk[\Z/2\Z]$ such that $A_n = \kk[\Z/2\Z]\#_\sigma H_n$.            
        \end{prop}
        \begin{proof}
            Consider the map $\psi$ defined in \cref{prop:A_n_n_surjects_onto_H_n_n}. A direct calculation shows that $A_n^{co~ \psi} = \prescr{co~ \psi}{{A_n}} = \mathrm{span}\{1,a^n\}$. Hence $\psi$ is normal.             
        \end{proof}
    \begin{example}
        Consider the map $\psi:A_2\to H_2$. Then $A_2 = \kk[\Z/2\Z]\#_\sigma H_2$ where the $2$-cocycle $\sigma$ satisfies:
        $$
            \sigma(g,g) = a^2,\quad \sigma(x,x) = a^2-1,\quad \sigma(g,x)=\sigma(x,g)=0
        $$
        on the generators $g,x$ of $H_2$.
    \end{example}
    \subsection{The category \texorpdfstring{$A_n-mod$}{An-mod}}
    Let $\mu$ be a primitive $2n$-th root of unity such that $\mu^2 = \omega$ and consider the Hopf algebra $A_n$. It is easy to see that as a vector space, $A_n$ has a basis $$\{\,a^iz^j\,|\,0\le i\le 2n-1,\, 0\le j\le n-1\,\}.$$ 
 Next, the algebra $A_n$ 
        has $2n$ orthogonal primitive idempotents, given by 
        $$
            e_i := \frac{1}{2n}\sum_{j=1}^{2n}\mu^{-ij}a^j
        $$
        with the action $a.e_i = \mu^i e_i$. As a left module over itself, we have a decomposition of $A_n$ into $A_n\cong \bigoplus_i A_ne_i$ where  $$P_i:=A_ne_i = span\{z^je_i\}_{j=0}^{n-1}$$ is an indecomposable projective module. The action of $z$ on $z^{n-1}e_i$ is 
        $$
            z. (z^{n-1}e_i) = z^n. e_i = (a^n-1)e_i = (\mu^{in}-1)e_i = \begin{cases}
            0 & \text{ if }\, i\equiv 0\mod 2\\ -2e_i & \text{ if }\, i\equiv 1\mod 2
            \end{cases}
        $$
        Let $i\in 2\Z$. Then $L_{i} := \operatorname{cosoc}(P_i)$ is a  simple $1$-dimensional $A_{n}$-module. The socle of $P_i$ is $L_{i-(2(n-1))}$. If $i\ne{j}$ are even, then $L_i\not\cong L_j$ as $A_{n}$-modules, because the action of $a$ on $L_i$ is given by multiplication by $\mu^i$. 

        \begin{lemma}\label{lem: P_i is simple for odd i}
            For $i \in 2\Z+1$, the projective $A_{n}$-module $P_i:=A_{n}e_i$ is simple.  
        \end{lemma} 
        \begin{proof}
            The element $a$ acts on $P_i$ by a diagonalizable operator with $n$ distinct eigenvalues, the eigenspaces being $\kk z^{j}e_i$ for $0\leq j\leq n-1$. Thus any $A_{n}$-submodule of $P_i$ is a direct sum of several of these eigenspaces. Since $z$ permutes these eigenspaces transitively, any non-trivial submodule of $P_i$ is the sum of all these eigenspaces. Hence $P_i$ is simple.
        \end{proof}
    
       Take $0\leq i\leq k\leq n$. Define the map\footnote{The indices of $P$ are taken modulo $2n$.} $f_{i, k}:P_{2i+1}\to P_{2k+1}$ by setting $f_{i, k}(z^te_{2i+1}) := z^{t+(k-i)}e_{2k+1}$ for $0\le{t}\le{n-1}$; in particular, this means that $f_{i, k}(z^{n-(k-i)}e_{2i+1})=-2e_{2k+1}$. 

        A direct calculation shows that this is indeed an $A_{n}$-module bijection, hence $A_n-mod$ has a single isomorphism class of non-invertible simple object $Q$, which is projective.

        To sum up, the isomorphism classes of simple $A_{n}$-modules are:
        \begin{itemize}
            \item $1$-dimensional modules $L_{2i}$ for $i=0,1, \ldots, n-1$ (these representations factor through the quotient map $A_n\twoheadrightarrow H_n$),
            \item a simple projective object $Q:=P_{1}$.
        \end{itemize}
    \begin{lemma}
        We have: $Q\otimes Q\cong\bigoplus_{i=0}^{n-1} P_{2i}$, so $[Q\otimes Q:Q]=0 $.
    \end{lemma}
    \begin{proof}
      From \cref{lemma: tensor product of simple projective object}(4) we see that $Q\otimes Q$ has $\bigoplus_{i=0}^{n-1} P_{2i}$ as a direct summand. Comparing the dimensions, we conclude that $Q\otimes Q\cong\bigoplus_{i=0}^{n-1} P_{2i}$.
    \end{proof}

        \begin{thm} \label{thm: basic projective extensions of taft algebras}
        Let $n\geq 2$. The category $A_n-mod$ with its forgetful functor is a FCA pair. Its Picard group is $\Z/n\Z$ and its extension constant is $r=0$.
        \end{thm}
  \begin{proof}
      We have determined that $A_n-mod$ is an admissible category with Picard group $\Z/n\Z$ and extension constant $r=0$. For this extension constant, we know (see Remark \ref{rmk:flip_compatible}) that every forgetful functor is flip-compatible, which proves the statement.
  \end{proof}
\section{First results on classification of FCA categories}\label{sec: classification first results}
    The full classification of almost-semisimple near-group categories (or even those admitting a fiber functor) is clearly a difficult question - for instance, it relies on the classification of pointed tensor categories. In this section, we take a modest first step of classifying FCA categories in which $\FPdim(Q)$ is prime.

    Throughout this section, $\D$ is an almost semisimple generalized near-group category admitting a fiber functor $F:\D\to \Vect$, with generalized fusion rule $(G,r)$. We denote by $Q$ the unique (up to isomorphism) non-invertible simple object with $q:=\FPdim(Q)$. We denote by $\C$ the full Serre subcategory generated by the invertible objects and $G$ stands for the Picard group $\pic(\D)$. Let $H:=\End(F)$ and $K:=\End(F\rvert_\C)$ be the corresponding Hopf algebras.
    \begin{prop}\label{prop: almost-semisimple near-group categories with prime dimension}
        If $q$ is prime, then one of the following statements holds:
    \begin{enumerate}
        \item $\D$ is a fusion near-group category with fusion rule $(\Z/q\Z,q-1)$ and  $K\cong \kk[\Z/q\Z]$.  
        \item $\D$ is a fusion near-group category with extension constant\footnote{Such near-group categories are known as Tambara-Yamagami near-group categories.} $r=0$ and $K \cong \kk[\Z/q\Z\times \Z/q\Z]$. 
      
        \item $\D$ is non-semisimple, $K$ is the Taft algebra of order $q$, and $r=0$. In this case, $H$ is a central extension of $K$ by $\kk[\Z/2\Z]$, as in \cref{thm: integral almost semisimple with r=0}.
        \end{enumerate}
    \end{prop}
    \begin{proof}
   
        Recall that by \cref{cor:divisibility} we have: $q\rvert (nd)$ and $(nd) \rvert q^2$. By \cref{eq:q_square}, we have: $r=q-\frac{nd}{q}$.
        
        Furthermore, if $n=1$ then $d=1$ since $Ext^1(\1, \1)=0$ (see \cite[Theorem 4.4.1]{EGNO}) and $Ext^1(\1, Q)=0$. 

        Assume that $q$ is prime. We get $$(n,d)\in \{(q^2, 1), (q,q), (q,1)\}.$$

        Let us first consider the case $d=1$. In this case $\1$ is projective, hence $\D$ is semisimple - that is, a fusion near-group category. In this case, $K$ is a semisimple Hopf algebra of dimension $n$. 
\begin{itemize}
    \item  If $n=q$ then $r=q-1$ and $K$ is a semisimple Hopf algebra of prime dimension, so $K\cong \kk[\Z/q\Z]$. 
    \item If $n=q^2$ then $r=0$ (so $\D$ is a Tambara-Yamagami near-group category). In this case, $\dim(K)=q^2$ so $K$ is isomorphic to either $\kk[\Z/q^2\Z]$ or to $\kk[\Z/q\Z\times \Z/q\Z]$. If $K\cong \kk[\Z/q^2\Z]$ then $G=\Z/q^2\Z$. In this case, $H^2(G, \kk^\times)=0$ and \cref{prop: dim(Q) in non-semisimple integral near-group categories} implies that $|G|=q^2 | q$, which is a contradiction. This rules out $K = \kk[\Z/q^2\Z]$.
  
\end{itemize}

      Finally, we consider the case $n=d=q$ (so $r=0$). In this case $\C$ (and $\D$) are non-semisimple, and $K$ is a non-semisimple Hopf algebra of dimension $q^2$. By \cite{Ng} it is isomorphic to the Taft algebra of order $n=q$. The remaining statement follows from \cref{thm: integral almost semisimple with r=0}, and the fact that for $r=0$ every fiber functor is flip-compatible (see Remark \ref{rmk:flip_compatible}).

      
    
    

    \end{proof}  
    All three types of categories listed in \cref{prop: almost-semisimple near-group categories with prime dimension} may occur:
    \begin{example}
    \mbox{}
        \begin{enumerate}
            \item The category $\Rep(S_3)$ is an example of type (1). Here $S_3$ stands for the symmetric group on $3$ symbols.
            \item The category $\Rep(D_8)$ is an example of type (2). Here $D_8$ is the dihedral group of order $8$. 
            \item The category $\Rep(A_p)$ is an example of type (3). Here $A_p$ is the extension of a Taft algebra of order $p$, as defined in \cref{sec: Extending H_{n,d}}.
        \end{enumerate}
    \end{example}

\begin{remark}
    There exists only one non-semisimple Hopf algebra $H$, up to isomorphism, such that $(\Rep(H), Forget)$ is an FCA pair with $\FPdim(Q) =q = 2$, by \cref{prop: almost-semisimple near-group categories with prime dimension} and \cite{Ste}. This is the Hopf algebra\footnote{In \cite{Ste}, this Hopf algebra is denoted by $A''_{C_4}$.} $A_2$ from \cref{sec: Extending H_{n,d}}.
\end{remark}
\printbibliography

\end{document}